\documentclass{siamltex}

\pdfoutput = 1
\usepackage[T1]{fontenc}
\usepackage[latin1]{inputenc}
\usepackage{latexsym}
\usepackage[pdftex]{graphicx}
\usepackage{stmaryrd}
\usepackage{amsmath,amssymb}
\usepackage[activate={true,nocompatibility},spacing,kerning]{microtype}
\usepackage{remreset}
\usepackage[pdftex,
            colorlinks,
            hyperfootnotes=false,
            pdffitwindow=true,
            plainpages=false,
            pdfpagelabels=true,
            pdfpagemode=UseOutlines,
            pdfpagelayout=SinglePage,
            pdftitle={Asymptotic Independence of the Extreme Eigenvalues of GUE},
            pdfauthor={Folkmar Bornemann, Technische Universität München},
            pagebackref,
            hyperindex,
            filecolor=purple,
            urlcolor=purple,
            citecolor=blue,
            linkcolor=blue]{hyperref}
\usepackage{myharvard}
\usepackage{bm}
\usepackage[sc,osf]{mathpazo}

\makeatletter
\@removefromreset{figure}{section}
\@removefromreset{equation}{section}
\makeatother

\input{macros.sty}

\title{Asymptotic Independence of the Extreme Eigenvalues of GUE}

\author{Folkmar Bornemann\thanks{Zentrum Mathematik, Technische Universität München,
        Boltzmannstr. 3, 85747 Garching, Germany ({\tt bornemann@ma.tum.de}). Manuscript
        as of \today.}}

\begin{document}

\maketitle

\begin{abstract} We give a short, operator-theoretic proof of the asymptotic independence (including a first correction term) of the
minimal and maximal eigenvalue of the $n \times n$ Gaussian Unitary Ensemble in the large matrix limit $n \to \infty$. This is done by representing
the joint probability distribution of the extreme eigenvalues as the Fredholm determinant of an operator matrix that
asymptotically becomes diagonal. As
a corollary we obtain that the correlation of the extreme eigenvalues asymptotically behaves like~$n^{-2/3}/4\sigma^2$,
where $\sigma^2$~denotes the variance of the Tracy--Widom distribution. While we conjecture that the extreme eigenvalues
are asymptotically independent for Wigner random hermitian matrix ensembles in general, the actual constant in the asymptotic behavior of the correlation
turns out to be specific and can thus be used to  distinguish the Gaussian Unitary Ensemble statistically
from other Wigner ensembles.
\end{abstract}

\section{Introduction}
We consider the $n\times n$ Gaussian Unitary Ensemble (GUE) with the joint probability distribution of its (unordered) eigenvalues given by
\[
p_n(\lambda_1,\ldots,\lambda_n) = c_n e^{-\lambda_1^2-\cdots-\lambda_n^2} \prod_{i<j} |\lambda_i-\lambda_j|^2
\]~\\*[-3mm]
and denote the induced minimal and maximal eigenvalue by $\lambda_{\min}^{(n)}$ and $\lambda_{\max}^{(n)}$.
\citeasnoun{Najim} have recently shown
the asymptotic independence of the edge-scaled extreme eigenvalues, that is, they proved
\begin{equation}\label{eq:najim}
\prob\left(\tilde\lambda_{\min}^{(n)} \leq x,\; \tilde\lambda_{\max}^{(n)} \leq y\right) =
\prob\left(\tilde\lambda_{\min}^{(n)} \leq x\right) \cdot \prob\left( \tilde\lambda_{\max}^{(n)} \leq y\right) + o(1)\qquad (n\to\infty)
\end{equation}
with the fluctuations
\[
\tilde\lambda_{\min}^{(n)} = 2^{1/2} n^{1/6} \left(\lambda_{\min}^{(n)} + \sqrt{2n}\right),\quad  \tilde\lambda_{\max}^{(n)} = 2^{1/2} n^{1/6} \left(\lambda_{\max}^{(n)} - \sqrt{2n}\right).
\]
The asymptotic independence can been used \cite{Najim2} to design, based on the \emph{ratio} of the extreme eigenvalues, a statistical test for the randomness of matrices that does \emph{not}
depend on estimating the actual variance of the distribution of the matrix entries (that is, the unknown level of noise in some applications).

In this paper we shall improve upon these results by showing that the \emph{correlation} of the extreme eigenvalues is a simple, scale-independent device to distinguish the GUE
statistically from other Wigner random hermitian matrix ensembles (and not just from non-random matrices like the ratio-based test).
To this end we establish a first correction term to the asymptotic independence
(\ref{eq:najim}), namely
\begin{multline}\label{eq:result}
\prob\left(\tilde\lambda_{\min}^{(n)} \leq x,\; \tilde\lambda_{\max}^{(n)} \leq y\right) \\
= \prob\left(\tilde\lambda_{\min}^{(n)} \leq x\right) \cdot \prob\left( \tilde\lambda_{\max}^{(n)} \leq y\right) +\tfrac14 F_2'(-x) F_2'(y) n^{-2/3} + O(n^{-4/3})
\end{multline}
as $n \to \infty$, locally uniform in $x$ and $y$. Here, $F_2$ denotes the Tracy--Widom distribution, see (\ref{eq:TW}) below. In fact,
the correction term comes as an additional benefit from a short and conceptually simple new proof of the asymptotic independence that  explains it straightforwardly from the
asymptotic diagonalization of a certain operator
matrix. In contrast, \citename{Najim} based their original proof on quite a detailed and lengthy study of
the classical power series of the Fredholm determinants representing the probability distributions in (\ref{eq:najim}).

In Section~\ref{sect:corr} we discuss the correlation of the extreme eigenvalues. In Section~\ref{sect:universality} we comment on universality.
Finally, in Section~\ref{sect:proof}, we prove the expansion (\ref{eq:result}).

\section{The Correlation of the Extreme Eigenvalues of GUE}\label{sect:corr}
\begin{figure}[tbp]
\begin{center}
\begin{minipage}{0.8125\textwidth}
{\includegraphics[width=\textwidth]{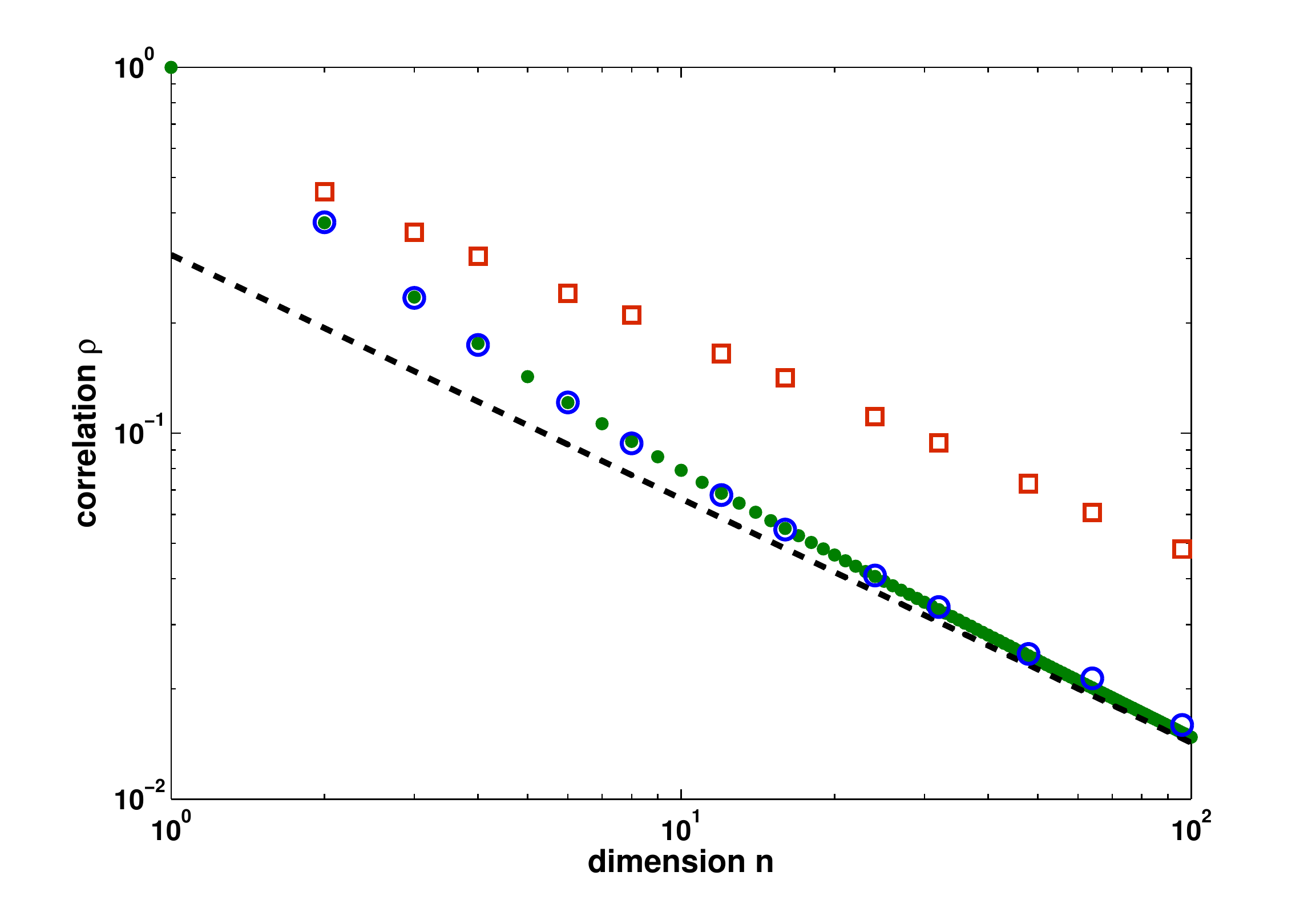}}
\end{minipage}
\end{center}\vspace*{-0.0625cm}
\caption{The dots show the values of the correlation $\rho$ of the extreme eigenvalues
of the $n\times n$ GUE as obtained from a numerical evaluation of the Fredholm determinant (\ref{eq:det}) by the method of \protect\citeasnoun{Bornemann}.
The dashed line shows the leading order term $n^{-2/3}/4\sigma^2$ of the asymptotic expansion (\ref{eq:rhoex}). The circles show the sample correlation
for $10^6$ realizations of $n\times n$ matrices drawn from the GUE. To compare with, the squares show the same
for $10^6$ realizations of $n\times n$ hermitian matrices whose algebraic degrees of freedom are uniformly distributed on $[-1,1]$.}
\label{fig:rho}
\end{figure}

Since both, $F_2'(-x)$ and $F_2'(y)$, are probability densities it follows from (\ref{eq:result}) that the covariance of the edge-scaled extreme eigenvalues of the GUE satisfies
\[
\cov\left(\tilde\lambda_{\min}^{(n)},\,\tilde\lambda_{\max}^{(n)}\right) = \tfrac14 n^{-2/3} + O(n^{-4/3}) \qquad (n\to \infty).
\]
Therefore, because of scale and shift invariance and by recalling (\ref{eq:edgeworth}) below, we get the correlation of the \emph{unscaled} extreme
eigenvalues (or of \emph{any} rescaling thereof) as
\begin{equation}\label{eq:rhoex}
\rho\left(\lambda_{\min}^{(n)},\,\lambda_{\max}^{(n)}\right)  = \frac{n^{-2/3}}{4\sigma^2}  + O(n^{-4/3}) \qquad (n\to \infty),
\end{equation}
where $\sigma^2 = 0.81319\,47928\,32957\cdots$ is the variance of the Tracy--Widom distribution. Figure~\ref{fig:rho} visualizes that the leading order term of
this expansion is actually quite a precise approximation of the correlation even for rather small dimensions $n$.

We observe a different asymptotic behavior for random hermitian matrices whose algebraic degrees of freedom are \emph{uniformly} distributed on $[-1,1]$.
Though the data shown in Figure~\ref{fig:rho} hint at an asymptotic behavior of the correlation of the form $\rho \simeq c\, n^{-2/3}$ here too, the constant $c$ is now,
quite distinguishably, about three times as large as for the GUE. Therefore, the correlation of the extreme eigenvalues may be used as a simple and effective
scale-independent device to distinguish the GUE statistically from other Wigner ensembles \citeaffixed{MR1727234}{as defined in}.

\section{Universality}\label{sect:universality}

Within the class of Wigner random hermitian matrix ensembles there are several limit laws known to hold {\em universally}.
Examples are the universality of the limit eigenvalue density, as given by Wigner's semicircle law, and of the limit distribution of the
(properly rescaled) fluctuations of the maximal eigenvalue, as given by the Tracy--Widom distribution \citeaffixed{MR1727234}{see}.
It is therefore reasonable to \emph{conjecture} the universality of the asymptotic independence of the extreme eigenvalues. In fact, the
sample correlation (squares) shown in Figure~\ref{fig:rho} for a concrete non-Gaussian Wigner ensemble strongly points into that direction.

However, since the asymptotic behavior $\rho \simeq c\, n^{-2/3}$ observed for this example differs in the constant $c$, it appears that the correction term in (\ref{eq:result}) has to be specific to the GUE.
We offer the following explanation for this effect. Choup's \citeyear{MR2233711,MR2406805} Edgeworth expansion for the GUE, that is,
\begin{equation}\label{eq:edgeworth}
\prob\left( \tilde\lambda_{\max}^{(n)} \leq t\right) = F_2(t)  + \gamma(t) n^{-2/3} + O(n^{-1})
\end{equation}
\citeaffixed[Thm.~1.3]{MR2406805}{where the coefficient $\gamma(t)$ is actually given by an explicit, though quite lengthy expression, see}, allows us to infer from (\ref{eq:result}) a likewise Edgeworth expansion of the joint probability distribution, namely
\begin{multline}\label{eq:jointedgeworth}
\prob\left(\tilde\lambda_{\min}^{(n)} \leq x,\; \tilde\lambda_{\max}^{(n)} \leq y\right)
= \left(1-F_2(-x)\right)\cdot F_2(y)  \\*[2mm] + \left((1-F_2(-x))\gamma(y)-\gamma(-x) F_2(y) + \tfrac14 F_2'(-x)F_2'(y)\right) n^{-2/3} + O(n^{-1}).
\end{multline}
(Note the considerable amount of cancelation that would have taken place within the order $O(n^{-2/3})$ terms if we had established (\ref{eq:result}) from those Edgeworth expansions at the first hand.)
Though the leading order term $F_2(t)$ of the Edgeworth expansion~(\ref{eq:edgeworth}) is known to be universal, the coefficient $\gamma(t)$ of the first correction
will  in general, as in the central limit theorem,  depend on some higher order moments of the underlying distribution of the matrix entries. Now, since the correction term to the asymptotic independence in (\ref{eq:result})
is contributing to exactly the same level of approximation in the expansion (\ref{eq:jointedgeworth}), namely to the order $O(n^{-2/3})$ term,
it will also most likely in general depend on the specific probability distribution of the matrix entries.

\section{Proof of the Asymptotic Expansion (\ref{eq:result})}\label{sect:proof}

Starting point is the well known representation
\[
p_n(\lambda_1,\ldots,\lambda_n) = \frac{1}{n!} \det( K_n(\lambda_i,\lambda_j))_{i,j=1}^n
\]
of the joint eigenvalue distribution in terms of the finite rank kernel (the second equality follows from the Christoffel--Darboux formula)
\begin{equation}\label{eq:Kn}
K_n(\xi,\eta) = \sum_{k=0}^{n-1} \phi_k(\xi)\phi_k(\eta) = \frac12\frac{\phi_n(\xi)\phi_n'(\eta) - \phi_n'(\xi)\phi_n(\eta)}{\xi-\eta} - \frac12 \phi_n(\xi)\phi_n(\eta)
\end{equation}
that is built from the $L^2(\R)$-orthonormal system of the Hermite functions
\[
\phi_m(t) = \frac{e^{-t^2/2} H_m(t)}{\pi^{1/4}\sqrt{m!} \,2^{m/2}}.
\]
From this representation we get the determinantal formulae \citeaffixed[§5.4]{MR1677884}{see}
\begin{align*}
\prob\left(X \leq \lambda_{\min}^{(n)}\right) = \det\left(I - K_n \projected{L^2(-\infty,X)} \right)&, \quad
\prob\left(\lambda_{\max}^{(n)}\leq Y \right) = \det\left(I - K_n \projected{L^2(Y,\infty)} \right),\\*[3mm]
\prob\left(X \leq \lambda_{\min}^{(n)},\; \lambda_{\max}^{(n)} \leq Y\right) &= \det\left(I - K_n \projected{L^2((-\infty,X) \cup (Y,\infty))} \right)
\end{align*}
with the natural constraint $X< Y$ (otherwise the last probability would be zero).
Whereas \citename{Najim} discuss exactly these determinants using Fredholm's power series, we refer to the fact that, for $X < Y$,
\[
L^2((-\infty,X) \cup (Y,\infty)) = L^2(-\infty,X) \oplus L^2(Y,\infty)
\]
which implies (for the equivalence of a single integral operator on a union of disjoint intervals with a system of integral operators see
\citeasnoun{MR1744872}, §VI.6.1, or \citeasnoun{Bornemann}, §8.1)
\begin{equation}\label{eq:det}
\prob\left(X \leq \lambda_{\min}^{(n)},\; \lambda_{\max}^{(n)} \leq Y\right) = \det\left(I -
\begin{pmatrix}
K_n & K_n \\*[1mm]
K_n & K_n
\end{pmatrix}\projected{L^2(-\infty,X) \oplus L^2(Y,\infty)}\right).
\end{equation}
Now, edge-scaling, that is, $X=-\sqrt{2n} -2^{-1/2}n^{-1/6} x$ and $Y=\sqrt{2n}+2^{-1/2}n^{-1/6}y$, transforms the kernel entries into
\begin{align*}
K_{11}^{(n)}(\xi,\eta) & = 2^{-1/2} n^{-1/6} K_n(-\sqrt{2n} -2^{-1/2}n^{-1/6} \xi, -\sqrt{2n} -2^{-1/2}n^{-1/6} \eta)\\*[1mm]
K_{12}^{(n)}(\xi,\eta) & = 2^{-1/2} n^{-1/6} K_n(-\sqrt{2n} -2^{-1/2}n^{-1/6} \xi, \sqrt{2n} +2^{-1/2}n^{-1/6} \eta)\\*[1mm]
K_{21}^{(n)}(\xi,\eta) & = 2^{-1/2} n^{-1/6} K_n(\sqrt{2n}+2^{-1/2}n^{-1/6} \xi, -\sqrt{2n} -2^{-1/2}n^{-1/6} \eta)\\*[1mm]
K_{22}^{(n)}(\xi,\eta) & = 2^{-1/2} n^{-1/6} K_n(\sqrt{2n} +2^{-1/2}n^{-1/6} \xi, \sqrt{2n} +2^{-1/2}n^{-1/6} \eta)
\end{align*}
and we obtain (note that $X<Y$ eventually for $n$ large, if $x$ and $y$ stay bounded)
\begin{multline}\label{eq:master}
\prob\left(-\tilde\lambda_{\min}^{(n)} \leq x\right) = \det\left(I - P_x\, K_{11}^{(n)} P_x\right),\quad  \prob\left(\tilde\lambda_{\max}^{(n)}\leq y \right) =
\det\left(I - P_y\, K_{22}^{(n)} P_y\right),\\*[4mm]
\shoveleft \qquad\quad \prob\left(-\tilde\lambda_{\min}^{(n)} \leq x,\; \tilde\lambda_{\max}^{(n)} \leq y\right) = \det\left(I -
\begin{pmatrix}
 P_x\, K_{11}^{(n)} P_x & P_x \,K_{12}^{(n)} P_y \\*[1mm]
P_y\, K_{21}^{(n)}P_x & P_y \,K_{22}^{(n)}P_y
\end{pmatrix}\right).\quad\;\,\end{multline}
Here, the operator matrix operates on the space $L^2(\R)\oplus L^2(\R)$ and the orthonormal projection $P_t : L^2(\R) \to L^2(t,\infty)$ is simply given by the multiplication operator
with the characteristic function $\chi_t$ of $(t,\infty)$. Plugging the Plancherel--Rotach
expansion \cite[Theorem~8.22.9]{MR0372517} of the Hermite functions, that is, the locally uniform expansion
\[
\phi_n\left(\sqrt{2n} + 2^{-1/2}n^{-1/6} t\right) = 2^{1/4}n^{-1/12} \left(\Ai(t) - \tfrac12\Ai'(t) n^{-1/3} + O(n^{-2/3})\right),
\]
into the rightmost expression defining the kernel $K_n$ in (\ref{eq:Kn}) yields the locally uniform asymptotic expansion
\[
K_{11}^{(n)}(\xi,\eta) = K_{22}^{(n)}(\xi,\eta) = K(\xi,\eta) + O(n^{-2/3})\qquad (n\to\infty)
\]
with the Airy kernel
\[
K(\xi,\eta) = \frac{\Ai(\xi)\Ai'(\eta) - \Ai'(\xi)\Ai(\eta)}{\xi-\eta}.
\]
Furthermore, \citeasnoun[Theorem~1.2]{MR2233711} proved that this expansion of kernels implies the asymptotic expansion
\[
P_t\, K_{11}^{(n)} P_t  = P_t\,K_{22}^{(n)}P_t = P_t K P_t + O(n^{-2/3})\qquad (n\to\infty)
\]
of the induced trace class operators (that is, the error $O(n^{-2/3})$ is also valid in trace norm). Completely analogously, by using the Plancherel--Rotach expansion once more
and recalling the symmetry $\phi_n(-t)=(-1)^n\phi_n(t)$, we obtain the locally uniform expansion
\[
K_{12}^{(n)}(\xi,\eta) = K_{21}^{(n)}(\xi,\eta) = \tfrac12(-1)^{n-1} \Ai(\xi)\Ai(\eta)\, n^{-1/3} + O(n^{-1})\qquad (n\to\infty),
\]
which, by the same arguments as \citename{MR2233711}'s, implies the validity of the asymptotic expansion
\[
P_t K_{12}^{(n)} P_s  = P_t K_{21}^{(n)} P_s = \tfrac12 (-1)^{n-1} (\chi_t \Ai \otimes \Ai \,\chi_s) \, n^{-1/3} + O(n^{-1})\qquad (n\to\infty)
\]
of the induced trace class operators. This shows that the off-diagonal operators in (\ref{eq:master}) have trace norm $O(n^{-1/3})$. Thus,
by the local Lipschitz continuity of the
determinant with respect to the trace norm \citeaffixed[Thm.~3.4]{MR2154153}{see} we get
\begin{multline*}
\prob\left(-\tilde\lambda_{\min}^{(n)} \leq x,\; \tilde\lambda_{\max}^{(n)} \leq y\right) = \det\left(I -
\begin{pmatrix}
 P_x\, K_{11}^{(n)} P_x & 0 \\*[1mm]
0 & P_y \,K_{22}^{(n)}P_y
\end{pmatrix} \right) + O(n^{-1/3})\\*[2mm]
= \det\left(I - P_x\, K_{11}^{(n)} P_x\right)\cdot  \det\left(I - P_y\, K_{22}^{(n)} P_y\right) + O(n^{-1/3})\\*[2mm]
= \prob\left(-\tilde\lambda_{\min}^{(n)} \leq x\right) \cdot \prob\left(\tilde\lambda_{\max}^{(n)}\leq y \right) + O(n^{-1/3}),
\end{multline*}
where we have used the multiplication rule of the determinant for \emph{diagonal} operator matrices.
This proves the asymptotic independence result (\ref{eq:najim}) of \citename{Najim}; here with an additional estimate of the order of approximation, however.

We now go one step
further and make the error term explicit to the leading order.
We start with the factorization
\begin{multline*}
\prob\left(-\tilde\lambda_{\min}^{(n)} \leq x,\; \tilde\lambda_{\max}^{(n)} \leq y\right) = \det\left(I -
\begin{pmatrix}
 P_x\, K_{11}^{(n)} P_x & 0 \\*[1mm]
0 & P_y \,K_{22}^{(n)}P_y
\end{pmatrix} \right) \cdot \\*[2mm]
\det\left(I -
\begin{pmatrix}
 0 & (I-P_x\, K_{11}^{(n)} P_x)^{-1}P_x \,K_{12}^{(n)} P_y  \\*[1mm]
(I-P_y\, K_{22}^{(n)} P_y)^{-1}P_y \,K_{21}^{(n)} P_x & 0
\end{pmatrix} \right).
\end{multline*}
As above, the first determinant evaluates to
\begin{equation}\label{eq:prod}
\det\left(I - P_x\, K_{11}^{(n)} P_x\right)\cdot  \det\left(I - P_y\, K_{22}^{(n)} P_y\right) = \prob\left(-\tilde\lambda_{\min}^{(n)} \leq x\right) \cdot \prob\left(\tilde\lambda_{\max}^{(n)}\leq y \right).
\end{equation}
The second determinant can be written for short as
\[
\det\left(I - \frac{(-1)^{n-1}n^{-1/3} }{2}
\begin{pmatrix}
 0 &   T_{xy} + O(n^{-2/3}) \\*[1mm]
  T_{yx}  + O(n^{-2/3}) & 0
\end{pmatrix} \right)
\]
with the rank-one operator
\[
T_{ts} = (I-P_t K P_t)^{-1} \chi_t \Ai \otimes \Ai\,\chi_s.
\]
By a straightforward operator decomposition \citeaffixed[(I.3.7)]{MR1744872}{see}, this evaluates and expands to
\begin{multline*}
\det\left(I- \tfrac14 n^{-2/3} \left( T_{xy} + O(n^{-2/3})\right)\left(T_{yx}  + O(n^{-2/3})\right)\right) \\*[2mm]
= \det\left(I- \tfrac14 n^{-2/3} T_{xy}T_{yx}\right) + O(n^{-4/3}).
\end{multline*}
The last determinant can actually be evaluated exactly. Indeed, by using the fact that $(f\otimes g)(v \otimes w) = \langle g,v\rangle f\otimes w$, and thus
\[
\det(I-(f\otimes g)(v \otimes w)) = \det(I - \langle g,v\rangle f\otimes w) = 1 - \langle g,v\rangle \langle f,w\rangle,
\]
we obtain
\[
\det\left(I- \tfrac14 n^{-2/3} T_{xy}T_{yx}\right)  = 1- \tfrac14 u(x)u(y)\, n^{-2/3}
\]
with the function
\begin{equation}\label{eq:u}
u(t) = \langle (I-P_t KP_t)^{-1} \chi_t \Ai,\Ai\,\chi_t\rangle.
\end{equation}
To summarize, our result so far is
\begin{multline*}
\prob\left(-\tilde\lambda_{\min}^{(n)} \leq x,\; \tilde\lambda_{\max}^{(n)} \leq y\right)\\*[2mm] =
\prob\left(-\tilde\lambda_{\min}^{(n)} \leq x\right) \cdot \prob\left(\tilde\lambda_{\max}^{(n)}\leq y \right)
\left(1-\tfrac14 u(x)u(y)\,n^{-2/3} + O(n^{-4/3})\right).
\end{multline*}
The product of the probabilities with the term $u(x) u(y)$ can further be simplified by expanding the first factor, written as the determinantal expression in (\ref{eq:prod}), through
\begin{equation}\label{eq:Edgeworth2}
\det\left(I-P_t\,K_{11}^{(n)}P_t \right)=\det\left(I-P_t\,K_{22}^{(n)}P_t \right) = \det(I-P_t KP_t) + O(n^{-2/3}).
\end{equation}
Now, by introducing the Tracy--Widom \citeyear{MR1257246} distribution
\begin{equation}\label{eq:TW}
F_2(t) = \det(I-P_t KP_t),
\end{equation}
and recalling that the function $u(t)$ defined in (\ref{eq:u}) actually satisfies $u(t) = F_2'(t)/F_2(t)$ \citeaffixed[p.~1132]{MR2054175}{see also}---a formula that was obtained in course
of Tracy and Widom's derivation of their famous Painlevé II representation of $F_2$---we finally get
\begin{multline*}
\prob\left(-\tilde\lambda_{\min}^{(n)} \leq x,\; \tilde\lambda_{\max}^{(n)} \leq y\right)\\*[2mm]
=\prob\left(-\tilde\lambda_{\min}^{(n)} \leq x\right) \cdot \prob\left(\tilde\lambda_{\max}^{(n)}\leq y \right)
- \tfrac14 F_2'(x)F_2'(y)\, n^{-2/3}  + O(n^{-4/3}),
\end{multline*}
which is easily seen to be equivalent to the asserted expansion (\ref{eq:result}).

\medskip

\paragraph{Remark} Numerical experiments using the methods of \citeasnoun{Bornemann} show that the error term of this expansion is indeed not better than of the order $O(n^{-4/3})$.

\subsection*{Acknowledgements} The author thanks Herbert Spohn for helpful comments on a first draft and for bringing up the issue of universality.

\bibliographystyle{kluwer}
\bibliography{article}

\end{document}